\documentclass{article}
\usepackage{amsmath}
\usepackage{amssymb}

\title{Subgroups of Clifford algebras}

\author{Robert Arnott Wilson}

\begin{document}

\maketitle

\begin{abstract}
Clifford algebras are used for constructing spin groups, and are therefore of particular importance in 
the theory of quantum mechanics. But the spin group is not the only subgroup of the Clifford algebra.
An algebraist's perspective on these groups and algebras may suggest ways in which they might be applied
 more widely  to describe the fundamental properties of matter.
 I do not claim to build a physical theory on top of the fundamental algebra, and my suggestions for
 possible physical interpretations are indicative only, and may not work.
 Nevertheless, both the existence of three generations of fermions and the symmetry-breaking of the weak interaction
 seem to emerge naturally from an extension of the Dirac algebra from complex numbers to quaternions.
\end{abstract}

\section{Introduction}
\subsection{Clifford algebras in the standard model}
The spin-type groups used in physics \cite{Griffiths,Zee,Woit} include not only the ordinary spin group
$Spin(3)\cong SU(2)$ and the relativistic \cite{Dirac} spin group $Spin(3,1)\cong SL(2,\mathbb C)$,
but also the gauge groups $Spin(2)\cong U(1)$ of quantum electrodynamics \cite{QED} and $Spin(3)\cong SU(2)$
of the weak interaction \cite{GWS}. Each of these groups arises from two different Clifford algebras \cite{Porteous},
so that already we have at least $6$ distinct Clifford algebras in the picture, namely $Cl(2,0)$, $Cl(0,2)$, $Cl(3,0)$,
$Cl(0,3)$, $Cl(3,1)$ and $Cl(1,3)$. 

The Dirac algebra, as used in the Feynman calculus in the standard model of particle physics, is a complex version of the real algebra
$Cl(1,3)$ generated by the Dirac gamma matrices $\gamma_1$, $\gamma_2$, $\gamma_3$ and $\gamma_0$, whose essential properties are
that they anti-commute with each other, and square to $-1$, $-1$, $-1$ and $+1$ respectively. Instead of using the scalar $i$ as a fifth generator,
it is conventional (or at least instructive) to use the matrix $\gamma_5:=i\gamma_0\gamma_1\gamma_2\gamma_3$, which anti-commutes with the original four gamma matrices and
squares to $+1$. With respect to these five gamma matrices as generators, therefore, the algebra has the structure of $Cl(2,3)$.

Another reasonable set of generators is
\begin{eqnarray}
\gamma_5,i\gamma_1,i\gamma_2,i\gamma_3,i\gamma_0
\end{eqnarray} 
which gives a structure $Cl(4,1)$, and a third possibility is
\begin{eqnarray}
\gamma_1,\gamma_2,\gamma_3,i\gamma_0,i\gamma_5
\end{eqnarray}
which gives a structure of $Cl(0,5)$. All these algebras are isomorphic (as algebras) to the algebra of all
$4\times 4$ complex matrices.
The same physical information can be packaged into slightly different mathematics by reversing the signature of the Clifford algebra.
We obtain $Cl(3,2)$, which consists of two copies of the algebra of $4\times 4$ real matrices, and $Cl(5,0)$ and $Cl(1,4)$, both of which consist
of two copies of $2\times 2$ quaternion matrices.

\subsection{Clifford algebras beyond the standard model}
Larger Clifford algebras have also been suggested \cite{Furey,perspective} for purposes such as including the strong force \cite{GM} gauge group $SU(3)$,
although $SU(3)$ is not itself a spin group. But $SU(3)$ is a subgroup of $SO(6)$, and by the Klein correspondence $Spin(6)\cong SU(4)$, which contains
$SU(3)$ as a subgroup. Hence this approach may also be useful in extending from the standard model to a `four-colour' model such as the
Pati--Salam model \cite{PatiSalam,Baez}. Similarly, Penrose's twistor theory \cite{twistors,twistorlectures} is based on $Spin(4,2)\cong SU(2,2)$, which can be constructed via the Clifford algebra $Cl(4,2)$ or $Cl(2,4)$. In fact, the algebras $Cl(3,3)$, $Cl(0,6)$ and $Cl(4,2)$ are all isomorphic to the algebra of $8\times 8$ real matrices, so all three approaches can be pursued simultaneously if desired. The reversed signatures $Cl(6,0)$ and $Cl(2,4)$ are both isomorphic to the
algebra of $4\times 4$ quaternion matrices, as are $Cl(5,1)$ and $Cl(1,5)$, so that these provide an alternative approach for those who appreciate the 
virtues of using quaternions \cite{Manogue}.

By now, I have mentioned 19 distinct Clifford algebras, and it is 
obvious that some clear mathematical principles are going to be necessary in order to make good choices for which ones to use in which physical contexts. The 
main principle I want to discuss in this paper
is the principle of choosing a particular Clifford algebra structure from the many that are in general available in a suitable abstract algebra. This is a \emph{mathematical} symmetry-breaking principle, which may or may not have anything to do with the
\emph{physical} process of `spontaneous symmetry-breaking'
that is usually invoked to explain the structure of the fundamental forces.

\section{Symmetry-breaking}
\subsection{Examples}
Every Clifford algebra is either a full matrix algebra over the real numbers, complex numbers or quaternions, or the direct sum of two such algebras,
isomorphic to each other. Therefore  the algebra contains much larger groups than the spin group it was designed to construct. Restricting from the large group of all invertible elements of the Clifford algebra, to the spin group, is a process of breaking a large symmetry group down to a small one.
For example, $Cl(0,2)$ is a quaternion algebra, so contains not only $Spin(2)\cong U(1)$, but also a copy of $SU(2)\cong Spin(3)$. 
Another example that I hope to show is relevant to physics is $Cl(3,1)$, that is a $4\times 4$ real matrix algebra, so contains not only $Spin(3,1)\cong
SL(2,\mathbb C)$ but also a copy of $SL(4,\mathbb R)\cong Spin(3,3)$.

In both these cases, the Clifford algebra breaks the symmetry. Thus $Cl(0,2)$ breaks the symmetry of $SU(2)$ by choosing a particular subgroup $U(1)$
to call $Spin(2)$. This example may be relevant to the symmetry-breaking of the weak gauge group $SU(2)$, that distinguishes between the neutral $Z$ boson and the two charged $W$ bosons.
Similarly, $Cl(3,1)$ breaks the symmetry of $SL(4,\mathbb R)$ by choosing a particular copy of $SL(2,\mathbb C)$ to call
$Spin(3,1)$. This example is not currently used in physics, as far as I am aware, but I hope to show that it is relevant to the breaking of
symmetry between the three generations of fermions. 

To see this operation of symmetry-breaking in a bigger context, we can use the Clifford algebra $Cl(3,3)$ to break the symmetry of the matrix group $SL(8,\mathbb R)$ down to $Spin(3,3)\cong SL(4,\mathbb R)$. Then we can use $Cl(3,1)$ to break $SL(4,\mathbb R)$ symmetry down to $Spin(3,1)\cong SL(2,\mathbb C)$.
The latter lies in both $Cl(3,0)$ and $Cl(1,2)$, which can be used to break the symmetry to either of the groups $Spin(3)\cong SU(2)$ or $Spin(1,2)\cong SL(2,\mathbb R)$, according to preference. 

If further symmetry-breaking is required, one can use $Cl(0,2)$ to break $SU(2)$ down to $U(1)$, and $Cl(2,0)$ to break $SL(2,\mathbb R)$ down to $U(1)$, or even use $Cl(1,1)$ to break $SL(2,\mathbb R)$ down to $GL(1,\mathbb R)$. These chains of symmetry-breaking
should be sufficient to deal with the main examples of symmetry-breaking in the standard model. The chain down through $Cl(3,1)$, $Cl(3,0)$ and $Cl(2,0)$ looks particularly promising.

\subsection{A quaternionic Dirac algebra}
One of many ways to construct $Cl(3,3)$ is to extend the complex numbers in the Dirac algebra to the quaternions. 
This is not just a matter of extending the $4\times 4$ complex matrices to $4\times 4$ quaternion matrices, but of adjoining 
quaternions $i,j,k$ that commute with
the 
four Dirac matrices $\gamma_1$, $\gamma_2$, $\gamma_3$ and $\gamma_0$. 
This process requires $8\times 8$ real matrices rather than $4\times 4$ quaternion matrices, since $j$ and $k$ 
must anti-commute with $\gamma_5$. To obtain anti-commuting generators required for a
Clifford algebra, we can take
\begin{eqnarray}
&i\gamma_1,i\gamma_2,i\gamma_3,i\gamma_0,j,k&
\end{eqnarray}
which exhibits the structure of $Cl(3,3)$. Some alternative sets of generators for $Cl(3,3)$ are
\begin{eqnarray}
&\gamma_1,\gamma_2,\gamma_3,\gamma_0,j\gamma_5,k\gamma_5;&\cr
&\gamma_1,\gamma_2,\gamma_3,\gamma_5,j\gamma_5,k\gamma_5.&
\end{eqnarray}
It is not immediately obvious which of these choices will be most useful for physical applications.
But the mathematical structure of $Cl(3,0)$ is quite different from that of $Cl(0,3)$, so the choice between $\gamma_1,\gamma_2,\gamma_3$ 
generating $Cl(0,3)$ and $i\gamma_1,i\gamma_2,i\gamma_3$ generating $Cl(3,0)$ is likely to be
significant.

As already noted, the signatures $(0,6)$ and $(4,2)$ are available in the same abstract algebra, for example
in the two sets of generators:
\begin{eqnarray}
&\gamma_1,\gamma_2,\gamma_3,i\gamma_0,j\gamma_0,k\gamma_0;&\cr
&i\gamma_1,i\gamma_2,i\gamma_3,\gamma_0,j\gamma_0\gamma_5,k\gamma_0\gamma_5.&
\end{eqnarray}
We therefore have a choice of three distinct spin groups
\begin{eqnarray}
Spin(3,3)&\cong& SL(4,\mathbb R)\cr
Spin(0,6)&\cong& SU(4)\cr
Spin(4,2)&\cong& SU(2,2)
\end{eqnarray}
but only the first of these is suitable for further symmetry-breaking with a smaller Clifford algebra. Both $Cl(3,1)$ and
$Cl(2,2)$ are available for this purpose, 
and although the former is the `obvious' one, the latter may also be useful.

The remaining spin group
\begin{eqnarray}
Spin(5,1)&\cong& SL(2,\mathbb H)
\end{eqnarray}
can similarly be used as input to further symmetry-breaking using any of $Cl(4,0)$, $Cl(1,3)$ or $Cl(0,4)$.
For the moment 
I will stick as closely as possible to the standard model, and work with $Cl(3,3)$ as an algebra of $8\times 8$
real matrices, with $Cl(3,1)$ used for symmetry-breaking of $Spin(3,3)$. I have suggested three possible generating sets,
two of which have $\gamma_5$ in the top (so-called pseudoscalar) degree, while the third has $\gamma_0$ in the top degree. The choice between
these options must be made on the basis of subgroups other than the spin group itself.

\section{Subgroups}

\subsection{Electroweak gauge groups}
In the standard model Dirac algebra $Cl(2,3)$, the modelling of electro-weak interactions uses the $2$-dimensional complex
subalgebra spanned by $1$ and $\gamma_5$. This permits the implementation of two copies of $U(1)$, generated as Lie groups
by $i$ and $i\gamma_5$. In other words, only a subgroup $U(1)$ of the weak gauge group $SU(2)$ appears in the Dirac algebra.
The full $SU(2)$ can only be implemented by extending to the quaternionic Dirac algebra $Cl(3,3)$ described above.

The Coleman--Mandula theorem \cite{Coleman} implies that the gauge group must lie in the centralizer algebra of the relativistic
spin group $SL(2,\mathbb C)$, generated by even products of the original four gamma matrices. This centralizer algebra is
generated by $j$, $k$ and $\gamma_5$, and is therefore isomorphic to $Cl(1,2)$, and thereby to the algebra of $2\times 2$
complex matrices. 

By choosing instead the generators $j\gamma_5$, $k\gamma_5$ and $\gamma_5$ we can identify this
algebra also as $Cl(3,0)$. 
Which generators should we choose? The four rotations are $i,j,k$ and $i\gamma_5$, and the three boosts are $\gamma_5$,
$j\gamma_5$ and $k\gamma_5$. Whatever happens, the symmetry is visibly broken, but we have a choice of which bits of symmetry
to keep, and which to abandon. This is largely a matter of physical interpretation, rather than mathematics, so it is not
appropriate to try to make such a decision here. What is clear, is that the four rotations give a group $U(1)$ generated by $i\gamma_5$,
and a group $SU(2)$ generated by $i,j,k$, and that these groups commute with each other, and generate a group $U(2)$. This is exactly the
formalism that  
is normally used to describe the unification of quantum electrodynamics (QED) with the weak force.

Thus consideration of the electro-weak interactions suggests that the third of the three suggested sets of generators for $Cl(3,3)$ is
likely to be the most useful.
There appear to be some minor technical differences between what we see here and what appears in the Feynman calculus, in that the roles
of $i$ and $i\gamma_5$ appear to have been interchanged. Since this can be adjusted by applying complex conjugation on half of the 
spinor, I suspect this is not a serious issue. Indeed, this difference may be more apparent than real:  it 
may well be possible to take the Feynman calculus for one generation exactly as it is,
and then to use $j\gamma_5$ or $k\gamma_5$ for the second and third generations, in place of $\gamma_5$.

There is complete symmetry amongst these three elements of the Clifford algebra, that is slightly obscured by the conventional notation.
If we write them instead as $i\gamma_0\gamma_1\gamma_2\gamma_3$,  $-k\gamma_0\gamma_1\gamma_2\gamma_3$ and $j\gamma_0\gamma_1\gamma_2\gamma_3$,  the symmetry becomes obvious.
There is also no clearcut answer to the question whether we should use compact $SU(2)$
or split $SL(2,\mathbb R)$ or complex $SL(2,\mathbb C)$ as the gauge group. This ambiguity also exists in the standard model, so again,
there should not be a problem with this. Indeed, the groups themselves play a rather minor role in the theory: it is the Clifford algebra
where all the calculations take place.

It is 
clear, I think, 
that electro-weak mixing can be expressed mathematically, 
not as a process of `mixing' groups together, or mixing Lie algebras together, but as a process of choosing 
bases for the matrix algebra to turn it into one or more Clifford algebras.
This viewpoint has 
some consequences for how to think about gauge groups in the general case. Although they may sometimes appear 
in the formalism as
groups or Lie algebras, it is certainly possible that these mathematical objects are secondary to 
a primary manifestation as a Clifford algebra.

\subsection{Spin groups}
The ordinary spin group $Spin(3)\cong SU(2)$ is generated as a Lie group by $\gamma_1\gamma_2$, $\gamma_2\gamma_3$ and
$\gamma_3\gamma_1$. It can be regarded either as the even part of $Cl(3,0)$, generated by $i\gamma_1,i\gamma_2,i\gamma_3$,
or as the even part of $Cl(0,3)$, generated by $\gamma_1,\gamma_2,\gamma_3$. The algebra $Cl(3,0)$ is isomorphic to the algebra
of all $2\times 2$ complex matrices, and therefore contains a subgroup $SL(2,\mathbb C)$ that is isomorphic to, but not equal to,
the relativistic spin group.  
I do not know if this copy of the group $SL(2,\mathbb C)$ has any reasonable physical interpretation.

The algebra $Cl(0,3)$, on the other hand, is the direct sum of two copies of the quaternion algebra. Generators for these two copies
may be taken as:
\begin{eqnarray}
\gamma_1+\gamma_2\gamma_3, \gamma_2+\gamma_3\gamma_1, \gamma_3+\gamma_1\gamma_2;\cr
\gamma_1+\gamma_3\gamma_2, \gamma_2+\gamma_1\gamma_3, \gamma_3+\gamma_2\gamma_1.
\end{eqnarray}
They therefore have a fairly obvious interpretation as left-handed and right-handed spins. Exponentiating these generators in the usual way
we obtain two copies of $SU(2)$, one left-handed and one right-handed. 
These two copies of $SU(2)$ intersect trivially, and therefore they generate a copy of $Spin(4)$. 

The same effect is achieved in the standard model using $Cl(4,0)$ generated by $\gamma_0,i\gamma_1,i\gamma_2,i\gamma_3$, with
pseudoscalar element $\gamma_5$. Here the full Clifford algebra is a $2\times 2$ quaternion matrix algebra, and the even part is the sum of
two quaternion algebras, which can be separated by the projections with the 
idempotents $(1\pm\gamma_5)/2$. The generators
can be taken as
\begin{eqnarray}
\gamma_1\gamma_2+i\gamma_0\gamma_3, \gamma_2\gamma_3+i\gamma_0\gamma_1, \gamma_3\gamma_1+i\gamma_2;\cr
\gamma_1\gamma_2-i\gamma_0\gamma_3, \gamma_2\gamma_3-i\gamma_0\gamma_1, \gamma_3\gamma_1-i\gamma_2.
\end{eqnarray}
The obvious question that occurs to an algebraist \cite{vdW} at this point is whether the extra complications in the standard model are
actually necessary, if the same effect can be achieved with the simpler mathematical machinery presented above? That is, can we take the same spin
terms as always, $\gamma_1\gamma_2$ and so on, but project with $(1\pm \gamma_1\gamma_2\gamma_3)/2$ without the
extra factor of $i\gamma_0$? I do not pretend to provide an answer to this, I merely ask the question.

The relativistic spin group $Spin(3,1)\cong SL(2,\mathbb C)$ can be obtained either from the even part of $Cl(3,1)$, generated
by $i\gamma_\mu$ for $\mu=0,1,2,3$, or from the even part of $Cl(1,3)$, generated by the $\gamma_\mu$ themselves. The algebra $Cl(3,1)$
is isomorphic to the $4\times 4$ real matrix algebra, while $Cl(1,3)$ is isomorphic to the $2\times 2$ quaternion matrix algebra,
and the even part in both cases is isomorphic to the $2\times 2$ complex matrix algebra. In both cases the pseudoscalar is
$\gamma_0\gamma_1\gamma_2\gamma_3$, which squares to $-1$. Hence there are no projections onto left-handed and
right-handed spins as such, but in the standard model this distinction is made instead by complex conjugation. This is done by 
extending the Clifford algebra either from $Cl(3,1)$ to $Cl(4,1)$, or from $Cl(1,3)$ to $Cl(2,3)$. Or, equivalently, complexifying the
Clifford algebra. 

So again, the question that
occurs to an algebraist is, is this extension really necessary, or is there already enough information in the odd part of $Cl(3,1)$ or $Cl(1,3)$?
In both cases, the Clifford algebra provides the necessary symmetry-breaking, either by choosing a copy of the complex numbers
inside the quaternions, or by imposing a complex $2$-space structure onto a real $4$-space. The discussion in this section suggests
that the algebra $Cl(1,3)$ provides a closer link to experimental properties of left-handed and right-handed spin, while the standard model
is somewhat ambivalent.

\subsection{The mass term in the Dirac equation} In the standard model Clifford 
algebra, whether it is regarded as $Cl(4,1)$ or $Cl(2,3)$
or $Cl(0,5)$, the pseudoscalar term is always $i$. This pseudoscalar is conventionally used for the mass term in the Dirac equation. 
But it may not be the most natural choice. The Dirac equation was modelled on a factorisation, into two linear factors, of Einstein's equation
\begin{eqnarray}
m^2c^2 &=& -p^2+E^2/c^2,
\end{eqnarray}
in that particular form, rather than, for example,
\begin{eqnarray}
E^2/c^2 &=& m^2c^2+p^2.
\end{eqnarray}
If the mass term is a scalar, then any re-arrangement of the equation introduces unphysical cross-terms in the squaring process, so is not possible. But the symmetrical form
\begin{eqnarray}
m^2c^2+p^2-E^2/c^2&=&0
\end{eqnarray}
suggests that we should really be looking at a group $SO(4,1)$, and therefore the Clifford algebra $Cl(4,1)$, with generators
$i\gamma_\mu$ for the energy and momentum terms, and therefore $\gamma_5$ for the mass term. This is consistent with the
earlier observation that extending to three generations of fermions seems to force the roles of $i$ and $i\gamma_5$ to be interchanged.

Indeed, if we want to extend to
a larger Clifford algebra in order to incorporate the three generations of fermions, then $i$ is no longer a natural choice for pseudoscalar.
The discussion so far has suggested either $\gamma_0$ or $\gamma_5$ as pseudoscalar, with perhaps a preference for $\gamma_0$.
The question really is whether (rest) mass ($\gamma_5$) or (total) energy ($\gamma_0$) is a better choice for pseudoscalar.
One has to make a choice on physical, not mathematical, grounds.
Experiment makes it clear that mass is not conserved in the weak interaction, while the whole of physics relies on the 
principle of conservation of energy. This suggests 
that taking $\gamma_0$ as the top degree may ultimately result in a more fundamental theory.

The current suggestion for generators for the quaternionic Dirac algebra is as follows, in two different notations, the second, more cumbersome, notation
 exhibiting the symmetry
between the three generations of fermions:
\begin{eqnarray}
&\gamma_1, \gamma_2,\gamma_3, &\gamma_5, j\gamma_5, k\gamma_5;\cr
&\gamma_1, \gamma_2,\gamma_3, &i\gamma_0\gamma_1\gamma_2\gamma_3, j\gamma_0\gamma_1\gamma_2\gamma_3, k\gamma_0\gamma_1\gamma_2\gamma_3.
\end{eqnarray}
This version has pseudoscalar $\gamma_0$, and therefore does not have an explicit energy term.
One can alternatively include the energy term as one of the generators, in order to obtain 
something closer to the standard model:
\begin{eqnarray}
\gamma_1, \gamma_2,\gamma_3, &\gamma_0, &j\gamma_5, k\gamma_5.
\end{eqnarray}
The latter breaks the $3+3$ symmetry to the $3+1+2$ symmetry that is a prominent feature of the standard model. But it has
pseudoscalar $\gamma_5$, and therefore treats rest mass as a more fundamental concept than energy.

At this point, we have a proposal for embedding the electro-weak gauge groups and the spin groups
into the quaternionic Dirac algebra $Cl(3,3)$, and the stage is set
for attempting to reproduce the rest of the standard model, with three generations of fermions, and
including the mixing angles and other parameters. 
Obtaining the values of the
parameters is not a matter of algebra, but obtaining a classification of the parameters is.

\section{The structure of the Clifford algebra}
\subsection{Undetermined parameters}
There are some parameters within the algebra $Cl(3,0)$, that describes electro-weak unification. 
These  parameters include one that describes the breaking of the symmetry within $SU(2)$, and three that describe how the pseudoscalar $U(1)$ mixes with the three dimensions of $SU(2)$. 

The remaining parameters describe how $Cl(3,0)$ relates to $Cl(0,3)$. 
The latter has $8$ dimensions in total, but under rotation symmetries of space
these reduce to $4$ different types, in degrees $0$, $1$, $2$ and $3$. The degree $0$ is the identity element of the algebra, so cannot mix with anything. Similarly, $Cl(3,0)$ has $8$ dimensions, and since the symmetry is broken, all $7$ non-identity elements participate 
individually in the mixing. Hence there are 
$3\times 7 = 21$ more parameters, 
making a total of $4+21$, 
exactly the number required for the standard model, according to some counts
at least. The $25$ parameters consist of $5\times 3$ mass terms, for 12 fermions and 3 bosons, $3+1$ mixing angles each from the
Cabibbo--Kobayashi--Maskawa (CKM) matrix \cite{Cabibbo,KM} and the Pontecorvo--Maki--Nakagawa--Sakata (PMNS) matrix \cite{Pontecorvo,MNS}, plus the fine-structure constant and the strong coupling constant.
  
  The main group of $21$ parameters splits under electroweak symmetries, according to my suggestion, into $3$ sets of $3+3+1$ of which one may contain lepton masses, and one quark masses, leaving the other to contain 
  $3+3+1$ of the $3+1+3+1$ mixing angles. The other two scalars here, I suggest, might be the fine-structure constant (to go with the lepton masses) and the strong coupling constant (to go with the quark masses). 
  That leaves the electroweak group of four parameters as the three boson masses (Higgs, $Z$ and $W$),
  together with the last mixing angle, which may be the CP-violating phase in the CKM matrix. This proposed allocation of parameters may not be exactly correct, of course, but it shows that the overall pattern of parameters matches closely to the standard model.
 
  In a little more detail, we can label the $3\times 7$ parameters with products of $\gamma_3$, $\gamma_1\gamma_2$, $\gamma_1\gamma_2\gamma_3$ with $i,j,k,\gamma_5,i\gamma_5,j\gamma_5,k\gamma_5$, and the remaining four with
  the products of $\gamma_5$ with $1,i,j,k$. 
  In the main block of $3\times 7$ parameters,
  the $i,j,k$ symmetry is a generation symmetry for fermions,
  while in the remaining four it represents the broken symmetry $1,k,j$ of the weak interaction, since there is a hidden multiplication by $i$
  in the theory.

  \subsection{A quaternionic notation}
  Since the notation inherited from Dirac obscures some of the symmetries, in particular the generation symmetry, I suggest an alternative notation.
  The algebra $Cl(1,3)$ is the algebra of $2\times 2$ quaternion matrices, so is the tensor product of quaternions with $2\times 2$ real matrices.
  Let us take $i',j',k'$ for these quaternions, to distinguish them from $i,j,k$ that have already been 
  used for a different quaternion algebra, and $2\times 2$ real matrices
  \begin{eqnarray}
  I:=\begin{pmatrix} 0&1\cr -1&0\end{pmatrix},&
  J:=\begin{pmatrix}0&1\cr 1&0\end{pmatrix},&
  K:=IJ=\begin{pmatrix}1&0\cr 0&-1\end{pmatrix}.
  \end{eqnarray}
  Then we can take Dirac matrices as generators, defined as follows:
  \begin{eqnarray}
  \gamma_1:=Ji', &\gamma_2:=Jj',&\gamma_3:=Jk',\quad \gamma_0:=K.
  \end{eqnarray}
  Our canonical copy of $Cl(0,3)$ is then generated by $Ji',Jj',Jk'$, and our canonical copy of $Cl(3,0)$ by $Ii, Ij, Ik$. Space symmetry is described by $i',j',k'$, and generation symmetry by $i,j,k$. 
  
  The pseudoscalar elements of these two subalgebras and the whole algebra are
  \begin{eqnarray}
  J&=&-\gamma_1\gamma_2\gamma_3,\cr
  I&=&i\gamma_5,\cr
  K&=&\gamma_0
  \end{eqnarray}
  respectively. One then sees the breaking of symmetry between the three generations, in the choice of $\gamma_5=-Ii$ in the standard model.
  
  The even part of the algebra corresponds to the diagonal matrices $1$ and $K$, while the odd part corresponds to the off-diagonal matrices $I$ and $J$.
  The ordinary spin group $SU(2)$ is generated by $i',j',k'$, and extends to the relativistic spin group generated by $Ii',Ij',Ik'$. The weak gauge
  group appears either as $SU(2)$ generated by $i,j,k$, or as $SL(2,\mathbb R)$ generated by $i,Ij,Ik$, 
  with a commuting copy of $U(1)$ generated by $I$.
  In particular, we have separated two concepts of complexification: one with $I$ and one with $i$. These two are mixed together in the standard model,
  in a way that depends on the choice of the first generation of fermions as special. 
  
    To describe the 25 parameters, we must choose an arbitrary direction in space, say the $k'$ direction, corresponding to the conventional
  choice of the $z$ direction in which to measure spin. Then we must mix the three elements $J, k', Jk'$ in $Cl(0,3)$ with the seven
  elements $i,j,k,I,Ii,Ij,Ik$ in $Cl(3,0)$. Hence the corresponding $21$ parameters lie in the elements
  \begin{eqnarray}
  \begin{array}{|ccc|ccc|c}
  Ji & Jj & Jk & Ki & Kj & Kk & K\cr
  k'i & k'j & k'k & Ik'i & Ik'j & Ik'k & Ik'\cr
  Jk'i & Jk'j & Jk'k & Kk'i & Kk'j & Kk'k & Kk'
  \end{array}
  \end{eqnarray}
  The remaining four parameters are perhaps best thought of as lying in $i, Ii, Ij, Ik$, and the whole system of $25$ parameters
  can then be arranged as follows:
   \begin{eqnarray}
  \begin{array}{|ccc|ccc|c}
  i&&&Ii&Ij&Ik&\cr\hline
  Ji & Jj & Jk & Ki & Kj & Kk & K\cr
  k'i & k'j & k'k & Ik'i & Ik'j & Ik'k & Ik'\cr
  Jk'i & Jk'j & Jk'k & Kk'i & Kk'j & Kk'k & Kk'
  \end{array}
  \end{eqnarray}
  
  The  second row of parameters does not have a momentum direction $k'$, and the first six are rotations, and contain generation labels
  $i,j,k$, so could reasonably be
  interpreted as the mixing angles between the three generations in the CKM and PMNS matrices. At a guess,
  $Ji, Jj, Jk$ lie in the CKM matrix, and $Ki, Kj, Kk$ in the PMNS matrix.   The first six elements in the 
  bottom two
  rows have momentum and generation labels, 
  and $9$ of them are boosts, so could reasonably be interpreted as
 $9$ of the $12$ fundamental fermion masses. The remaining $3$ are rotations, so perhaps represent the so far undetected neutrino masses.
  Of the remaining $7$ parameters, two are rotations, and I suggest $i$ to hold the CP-violating phase of the CKM matrix, and
  $Kk'$ to hold the CP-violating phase of the PMNS matrix, although it is perhaps equally plausible to allocate these the other way round. 
  
  This leaves five boosts, $K$, $Ik'$, $Ii$, $Ij$ and $Ik$, to hold the remaining
  three masses and two coupling constants. A reasonable guess might be that those involving $I$ are electroweak parameters, and $K$ is
  a strong force parameter, therefore the strong coupling constant. Perhaps $Ik'$ is the fine structure constant, and the remaining three
  the masses of the Higgs, $Z$ and $W$ bosons. But electroweak symmetry-breaking may significantly complicate this picture.
  
  \subsection{Even and odd}
  The even part of $Cl(3,3)$ consists of degrees $0,2,4$ and $6$, with dimensions $1$, $15$, $15$ and $1$ respectively. Bases are as follows,
  arranged with degrees $0$ and $4$ in the first block, and degrees $2$ and $6$ in the second:
  \begin{eqnarray}\begin{array}{c|ccc} 1& Ki&Kj&Kk\cr\hline Ki'&i'i&i'j&i'k\cr Kj'&j'i&j'j&j'k\cr Kk'&k'i&k'j&k'k\end{array} &&
   \begin{array}{c|ccc} K& i&j&k\cr\hline i'&Ki'i &Ki'j & Ki'k\cr j' & Kj'i&Kj'j&Kj'k\cr k'&Kk'i&Kk'j&Kk'k\end{array}
  \end{eqnarray}
  Neither block is closed under the Clifford multiplication, but the first block is closed under Jordan multiplication $AB+BA$, and the second
  is closed under Lie multiplication $AB-BA$. Indeed, the Lie bracket converts the degree $2$ part of the algebra into the Lie algebra of
  $Spin(3,3)$. Exponentiating the elements of degree $2$ therefore gives (the canonical copy of) $Spin(3,3)$. This is the only copy of
  $Spin(3,3)$ that preserves the Clifford algebra structure under its action by conjugation.
  
  The pseudoscalar $K$ defines two orthogonal central idempotents $(1\pm K)/2$ in the even part of the Clifford algebra, and the corresponding
  projections map onto two $16$-dimensional subalgebras, each of which is isomorphic to the full $4\times 4$ real matrix algebra. Since
 \begin{eqnarray}
 (1+K)(1-K) &=& 0
 \end{eqnarray}
  each subalgebra annihilates the other. In particular, each algebra contains a copy of $SL(4,\mathbb R)$ that acts on only half of the 
  $8$-dimensional real spinor. These subgroups may be useful for physics, but it must be stressed that they do not preserve the 
  Clifford algebra structure, and therefore they unavoidably break a great deal of symmetry. 
  They do, however, preserve the distinction between the even and odd parts of the Clifford algebra.
    
  The odd part of the algebra splits into degrees $1$, $3$ and $5$, with dimensions $6$, $20$ and $6$ respectively. 
  Degrees $1$ and $5$ have bases
  \begin{eqnarray}
 & Ji',Jj',Jk', Ii,Ij,Ik;&\cr
  &Ii',Ij',Ik',Ji,Jj,Jk&
  \end{eqnarray}
  respectively. As a representation of
$Spin(3,3)$, the degree $3$ part splits into a dual pair of $10$-dimensional real representations. These are obtained by multiplying the
$10$ elements
\begin{eqnarray}
\begin{array}{c|ccc}
1&&&\cr\hline
&i'i&i'j&i'k\cr
&j'i&j'j&j'k\cr
&k'i&k'j&k'k
\end{array}
\end{eqnarray}
by $J+I$ and $J-I$ respectively. This splitting corresponds to the splitting into the top-right and bottom-left $4\times 4$ blocks
of the $8\times 8$ matrices.
  
  Now
  \begin{eqnarray}
  (J\pm I)^2&=&0\cr
  (J+I)(J-I)&=& 2(1+K)\cr
  (J-I)(J+I)&=& 2(1-K)
  \end{eqnarray}
  so that each of these two $10$-spaces has Clifford product identically zero. 
  Their product in one order gives one projection of the even part of the algebra, and
  their product in the other order gives the other projection. 
  
  \subsection{A change of signature?}
  The unification of space and time in special relativity mixes momentum ($Ji'$,  $Jj'$, $Jk'$, 
  in the odd part of the Clifford algebra) with energy ($K$, in the even part).
  Similarly, the unification of the electroweak interactions mixes $i,j,k$ in the even part with $I$ in the odd part. In other words, both these unifications
  require groups that extend into the odd part of the Clifford algebra.
  This might mean that the Clifford algebra that I have suggested is irrelevant for physics. Or it might mean that the relationship between
  different Clifford algebra structures may throw some light on the process of symmetry-breaking.
  In particular, the necessity for using the odd part of the algebra for part of the Lorentz group, or relativistic spin group, and for part of the electro-weak gauge group, suggests the possibility that the strong gauge group $SU(3)$ may also lie partly in the odd part of the algebra.

Alternatively, we may prefer to change our choice of Clifford algebra structure to try to bring the standard model symmetries
into the even part of the algebra. Essentially, we have to make $I$ even and $K$ odd. One way to do this is to choose the following
generators, given in both notations:
\begin{eqnarray}
&\gamma_1,\gamma_2,\gamma_3,i\gamma_0, j\gamma_0,k\gamma_0;&\cr
&Ji',Jj',Jk',Ki,Kj,Kk.&
\end{eqnarray}
This has the effect of changing the signature to $(0,6)$, and the pseudoscalar to $I$ or $i\gamma_5$,
and changing the spin group to $Spin(6)\cong SU(4)$. 

There are now two obvious ways to collapse the Clifford algebra from a $6$-space to a $4$-space, by multiplying together
either the first three generators or the last three. Both give copies of $Cl(1,3)$. One of these has generators equal to the
original gamma matrices, and the corresponding spin group $Spin(1,3)$ is the usual relativistic spin group. The other one
has generators $i\gamma_0\gamma_5,i\gamma_0,j\gamma_0,k\gamma_0$, and the corresponding spin group is generated by
\begin{eqnarray}
i,j,k,\gamma_5,j\gamma_5,k\gamma_5.
\end{eqnarray}
Hence it describes electro-weak mixing in exactly the same way as I have already suggested. Both collapsed Clifford algebras
also have pseudoscalar $i\gamma_5$.

There is then an obvious copy of $U(3)$ inside the spin group, generated as a Lie group by the $9$ elements
\begin{eqnarray}
\begin{array}{ccc}
i+i', & Ii'i,& I(j'k+k'j),\cr
j+j', & Ij'j, & I(k'i+i'k),\cr
k+k', & Ik'k, & I(i'j+j'i).
\end{array}
\end{eqnarray}
The scalar copy of $U(1)$ is generated by $I(i'i+j'j+k'k)$, and $SU(3)$ is generated by the orthogonal complement of this.
There are other choices of $SU(3)$, but this is the most symmetrical.

Now comes the problem of interpreting the group $SU(3)$ constructed here. It is hard to support an interpretation as the 
gauge group of the strong interaction, 
since the mathematical properties of this copy of $SU(3)$ are not the same as those of the gauge group in the standard model.
In particular, it does not commute with the relativistic spin group, generated by
\begin{eqnarray}
&i',j',k',Ii',Ij',Ik',&
\end{eqnarray}
and thereby violates the Coleman--Mandula theorem. 
Moreover,
it does not commute with the weak gauge group generated by $i,j,k$. 

The only things it does commute with in the even part of
the algebra are two copies of the complex numbers, 
generated by $I$ and $I(i'i+j'j+k'k)$.
Nevertheless, this copy of $SU(3)$ may repay further study, so let us investigate its action on the Clifford algebra in more detail.

We have already seen the adjoint representation and a scalar in degree $2$, and the rest of the degree $2$ part 
is a $6$-dimensional real representation on 
the generators
\begin{eqnarray}
i'-i,j'-j,k'-k && I(j'k-k'j), I(k'i-i'k), I(i'j-j'i).
\end{eqnarray}
The latter can be given the structure of a complex $3$-space, using the natural scalar $I(i'i+j'j+k'k)$.
The degree $4$ part consists of equivalent pieces $8+1+6$, just
multiplied by the pseudoscalar $I$. In degrees $1$ and $5$ we again see copies of the $6$-dimensonal real representation,
and again there is a natural complex structure defined by $I(i'i+j'j+k'k)$. 

In degree $3$,
however, there are \emph{two} natural complex structures, defined by $I$ and $I(i'i+j'j+k'k)$.
Both give this $20$-dimensional space the structure of a complex
$10$-space, breaking into irreducibles as $1+3+6$. The complex scalars are contained in $J$ and $K$ times $(1+i'i+j'j+k'k)$.
The $3$-space consists of the anti-symmetric part
\begin{eqnarray}
&J(j'k-k'j), J(k'i-i'k), J(i'j-j'i),& \cr &K(j'k-k'j), K(k'i-i'k), K(i'j-j'i),&
\end{eqnarray}
and the $6$-space consists of the symmetric part
\begin{eqnarray}
&J(3-i'i-j'j-k'k), J(i'i-j'j), J(j'j-k'k), &\cr
& J(j'k+k'j), J(k'i+i'k), J(i'j+j'i),&\cr
&K(3-i'i-j'j-k'k), K(i'i-j'j), K(j'j-k'k), &\cr
& K(j'k+k'j), K(k'i+i'k), K(i'j+j'i).&
\end{eqnarray}

All of this structure is manifest in the quaternionic Dirac algebra, and therefore 
must surely have a sensible interpretation in physics.
But it seems quite clear that this copy of $SU(3)$ cannot be sensibly interpreted as the gauge group of quantum chromodynamics (QCD).
So what is it?

\section{Prospects for unification}

\subsection{Mass versus energy}
Mathematically, the distinction between $Cl(3,3)$ and $Cl(0,6)$ is most evident in the distinction between using $I$ or $K$ as
pseudoscalar, that is, in physics notation, the distinction between $i\gamma_5$ and $\gamma_0$. Physically this corresponds to 
a distinction between symmetry groups for which mass is scalar, and those for which energy is a scalar. Roughly speaking,
mass is a scalar in the theories of special relativity and (classical) electromagnetism, as well as quantum electrodynamics
and quantum chromodynamics, while energy is a scalar in general relativity, and in the theory of the weak interaction. 

The former (QED and QCD) are the quantum interactions for which the gauge bosons are massless.
The gauge groups $U(1)$ and $SU(3)$ therefore embed in the spin group $Spin(0,6)\cong SU(4)$ which is the subgroup
of the Clifford algebra that fixes the mass coordinate. 
But these groups do not necessarily commute with each other. The group $U(1)$ maps onto a subgroup $U(1)$ of $SU(3)$, 
or more generally of $U(3)$, and this
subgroup is defined by a $3\times 3$ unitary matrix, that appears in the standard model as the CKM matrix.

 In the weak interaction, on the other hand, the gauge bosons are massive, so that the gauge group
 does not act on mass as a scalar, but instead acts on energy as a scalar. 
 Hence the gauge group $SU(2)$ embeds instead in the spin group $Spin(3,3)\cong SL(4,\mathbb R)$ that fixes the
 energy coordinate.  The Pati--Salam model extends $SU(2)\cong Spin(3)$ acting on three coordinates to $Spin(4)\cong SU(2)\times SU(2)$
 acting on all six. The Clifford algebra formalism identifies the two scalars here, so that the group $Spin(4)$, which is not a subgroup
 of $SL(4,\mathbb R)$, is replaced by $SO(4)$, which is. In other words, the Pati--Salam model in effect uses both
 Clifford algebra structures $Cl(0,6)$ supporting $SU(4)$, and $Cl(3,3)$ supporting $SO(4)$ inside $SL(4,\mathbb R)$.
 But the model does not address the issue of symmetry-breaking, which arises from embedding both groups
 into $SL(8,\mathbb R)$, by identifying the two Clifford algebras as being two different ways of looking at the same thing.
 
 Similarly, in general relativity, in theory energy is a scalar (conserved, but not necessarily invariant), 
 so that the relevant group is $SL(4,\mathbb R)$,
 but most formulations assume that mass is a scalar (both invariant and conserved), 
which breaks the symmetry down to $SO(3,1)$. Now, as I have shown, in a flat spacetime it is not possible
 for both mass and energy to be scalars, and therefore general relativity incorporates a curvature of spacetime
 in order to solve this problem. On the other hand, if one abandons the assumption that mass is a scalar, which
 we must do in the case of the weak interaction anyway, then it may be possible to construct a version of general
 relativity in a flat spacetime, with gauge group $SL(4,\mathbb R)$, in such a way that the curvature of spacetime
 is replaced by variations in rest mass between different observers, accelerating with respect to each other.
 
In this way,  $Cl(3,3)$ and the associated group $Spin(3,3)\cong SL(4,\mathbb R)$
can potentially describe all the theory that one sees if one assumes that energy is a scalar, including the weak interaction, and
a theory of gravity that reduces to general relativity in the limit that rest masses can be treated as constants;
while $Cl(0,6)$ and the associated group $Spin(0,6)\cong SU(4)$ can describe all the theory one sees if one
assumes that mass is a scalar, including electrodynamics and the strong interaction.

Any unified theory that uses a Clifford algebra must reconcile this difference in viewpoint in some way, 
either by choosing one or the other of
mass and energy to be a scalar and allowing the other to vary, or by providing a translation between the two. 
This fact is already evident in the unification of QED (with scalar mass) and the weak force (with scalar energy). The translation
between the two is expressed in the standard model by a symmetry-breaking that arises from treating both mass and energy
as scalars simultaneously. 

In terms of both physical principles, and experimental reality, it surely makes more sense to treat energy as a scalar, rather
than mass. This is already done very effectively in the theory
of the weak interaction. But this 
can only be done by breaking the symmetry groups that rely on mass being
a scalar.
To do the same for the other forces
requires breaking the Lorentz symmetry group $SO(3,1)$ down to $SO(3)$,
and breaking $SU(3)$
down to essentially nothing using the CKM matrix.

\subsection{The strong gauge group}
If it is really possible to include all the unexplained parameters of the standard model inside the quaternionic Dirac algebra $Cl(3,3)$,
then this algebra must effectively contain the gauge group $SU(3)$ of the strong force as well. But there are some serious mathematical
problems with such a proposal, as this Clifford algebra is simply not big enough to contain a commuting product of all the gauge groups
and the relativistic spin group. 

The only copies of $SU(3)$ that are available break the symmetry significantly. But we have already seen,
in the case of electro-weak mixing, that the Clifford algebra imposes a symmetry-breaking that is not obvious in the standard model,
so maybe the same is true for the mixing of the strong force with the electroweak forces.
The triplet colour symmetry is unobservable, and the only available place to put it is in the triplet $\gamma_1\gamma_2, \gamma_2\gamma_3,
\gamma_3\gamma_1$ that is independent of the $25$ significant parameters.  This unavoidably identifies colour as being essentially the same
thing as the direction of spin.  The direction of a quark spin 
is surely not measurable, at least in practice, if not in theory.
The possible directions, moreover, form a $2$-parameter family, consistent with the property of colour confinement.

From a mathematical point of view, there is little prospect of 
linking the gauge group $SU(3)$ of quantum chromodynamics (QCD) as a theory of the strong force to a subgroup $SU(3)$ of the
Clifford algebra. The alternatives therefore seem to be either to keep the gauge group outside the Clifford algebra completely,
as in the standard model, or to try to re-write the formalism in terms of the $8$-dimensional subalgebra $Cl(0,3)$, generated by $\gamma_1,
\gamma_2,\gamma_3$, rather than the
$8$-dimensional group $SU(3)$. I have no idea whether this is possible, and therefore leave this as an open problem. It seems to require
a rather radical re-interpretation of the strong force as being not so much a force as a
quantisation of space
itself. Such an interpretation is certainly unconventional,  but is perhaps hinted at by the phenomenon of asymptotic freedom \cite{freedom1,freedom2}.

A more conventional approach would be to take the direct product of all the gauge groups and the spin group
\begin{eqnarray}
U(1) \times SU(2) \times SU(3) \times SL(2,\mathbb C)
\end{eqnarray}
and embed it in the real group
\begin{eqnarray}
SO(2) \times SO(4) \times SO(6) \times SL(4,\mathbb R)
\end{eqnarray}
so that there is an obvious $16$-dimensional real representation. Thus the group can be embedded in a Clifford algebra of an
$8$-dimensional space with signature $(8,0)$, $(5,3)$, $(4,4)$, $(1,7)$ or $(0,8)$. A breaking of symmetry into $2+6$ and $4+4$
can be accommodated in $Cl(4,3)$ or $Cl(0,7)$. Alternatively, a breaking of symmetry to incorporate the complex structure can be
accommodated in $Cl(7,0)$, $Cl(5,2)$, $Cl(3,4)$ or $Cl(1,6)$. 
All of these approaches, however, introduce a large number of
extra parameters whose meanings and values still have to be explained.

\subsection{Further remarks}
The suggested splitting as $3+1+2$ keeps $Cl(0,3)$ intact, but splits $Cl(3,0)$ into $Cl(1,0)$ and $Cl(2,0)$. A slightly closer match to the standard
model splitting into $SU(3)$, $U(1)$ and $SU(2)$ can be obtained by reversing the signatures, so that $Cl(0,3)$ splits into $Cl(0,1)$ and
$Cl(0,2)$, which are the complex numbers and quaternions, respectively. Hence we have obvious groups $U(1)$ and $SU(2)$, as subgroups
of the group $Spin(4)$ 
that lies inside $Cl(0,3)$. The other half of the algebra is then $Cl(3,0)$, the algebra of all $2\times 2$ complex
matrices, as a real $8$-dimensional algebra to replace the real $8$-dimensional group $SU(3)$. The group that we actually get is,
however, $GL(2,\mathbb C)$, also $8$-dimensional as a real Lie group, with a normal subgroup 
$SL(2,\mathbb C)$, that is isomorphic to, but distinct from, the relativistic spin group of Dirac. From a mathematical point of view,
I cannot see any strong argument for adopting this reversal of signature, but it is possible that are good physical arguments for it.

It is also interesting to consider the Georgi--Glashow model \cite{GG} in this context. They unified $3+2$ into $5$, and because they were working with unitary groups, they used $SU(5)$ for unification of $SU(2)$ and $SU(3)$. 
But this group has dimension $24$, and turned out to be too big, as it contains extra gauge bosons that cause proton decay, a phenomenon that has never been observed experimentally. If instead we unify Clifford algebras $Cl(2,0)$ and $Cl(0,3)$, then we obtain instead $Cl(2,3)$, that is
isomorphic to the Dirac algebra. The 
largest group that 
could possibly be relevant here is the 
general linear group $GL(4,\mathbb C)$ of (complex) dimension $16$. This group is  
small enough to
avoid the problem of proton decay, since it contains $13$ dimensions for the standard $13$ bosons, plus $3$ dimensions for
the spin group itself. 

A final remark on chirality is in order. The suggestion to use $i\gamma_5$ as a generator for the electromagnetic $U(1)$, rather than $i$,  
is affected by the fact that $i\gamma_5$ anti-commutes with $\gamma_\mu$ for $\mu=0,1,2,3$.
This implies that electromagnetism does not have a chirality, since the sign of $\gamma_1\gamma_2\gamma_3$ is not well-defined. On the other hand, the weak $SU(2)$ generated by $i,j,k$ commutes with these $\gamma_\mu$. This implies that the weak force does have a chirality, as is 
experimentally observed 
\cite{Wu}, since it distinguishes between $\gamma_1\gamma_2\gamma_3$ and its negative. Our choice of possible generators for the strong force algebra suggests
that the strong force also should not have a chirality, since the $\gamma_\mu$ anti-commute with each other.

  \subsection{Possible relationships with gravity}
At a philosophical level, the suggestion that the strong force is essentially describing the quantum structure of space relates closely to 
Einstein's approach \cite{GR1,GR2,GR3}
to gravity, namely that gravity is not so much a force as a description of the shape of spacetime on a macroscopic scale.
  At a mathematical level, the use of $Cl(3,3)$ and the spin group $Spin(3,3)\cong SL(4,\mathbb R)$ 
  may possibly be related to the use of
  $GL(4,\mathbb R)$ in some modern approaches to general relativity \cite{GL4R1,GL4R2}. But these superficial  
  observations are 
  a far cry from actually unifying the theories
  of general relativity and particle physics.
  
  A first remark is that in $Cl(3,3)$ there are many ways to find two copies of $SL(4,\mathbb R)$ that commute with each other. Therefore one could in principle use one
  copy of $SL(4,\mathbb R)$ for particle physics and the other for general relativity, simply placed side-by-side and not interacting with each other.
  Neither copy is equal to the spin group $Spin(3,3)$, and neither of them respects the Clifford algebra structure,
  but they are obtained   by suitable projections,
  defined by the pseudoscalar. 
  
    One could try to use the standard particle physics projections with $1\pm \gamma_5$, but the above discussion suggests that
  projections with $1\pm \gamma_0$ are likely to be more useful in this context. This reflects the fact that general relativity describes a gravity that
  depends on total energy ($\gamma_0$) rather than on rest mass ($\gamma_5$). If one tries to use both pairs of projections together, one 
  runs into problems because $\gamma_0$ and $\gamma_5$ do not commute with each other.
  Whatever choice of projections we make, the splitting of $Cl(3,3)$ into two copies of $4\times 4$ real matrices is a symmetry-breaking
  of the kind I have been discussing, down to $Cl(3,2)$. Further symmetry-breaking down to a single copy of $4\times 4$ real matrices
  can be achieved 
  by going to $Cl(3,1)$ or 
  $Cl(2,2)$.
  
  I have no idea whether it is actually possible to incorporate general relativity into $Cl(3,3)$ in this way. If it is, then the algebra implies
  that there is a certain overlap between the parts of the algebra that are used for particle physics and those that are used for gravity.
  Such an overlap is not part of current mainstream thinking on fundamental physics, but 
  has been considered at various times by Einstein \cite{Einsteininertia}, Sachs \cite{Sachs}, Penrose \cite{roadtoreality} and others \cite{uses}.

An alternative approach, more in line with mainstream ideas, is to go to a bigger Clifford algebra, so that this overlap can be avoided.
For example, $Cl(4,3)$ is the sum of two copies of $8\times 8$ real matrices, which permits one copy to be used for $Cl(3,3)$ to model
particle physics, and the other to act on four dimensions of spacetime and four dimensions of energy-momentum, to describe
classical physics and general relativity. Possibly it might be useful to go further, to $Cl(4,4)$, in order to have enough room for the
standard model strong gauge group $SU(3)$, rather than reducing this to a smaller group as I have suggested.

  \section{Conclusion}
  In this paper I have explored the subgroups of many small Clifford algebras in sufficient detail to conclude that $Cl(3,3)$ is the closest match among the Clifford algebras to the algebraic foundations of the electro-weak part of 
   the standard model of particle physics. I have shown that it incorporates the three generations of fermions
  in the minimal possible way, consistent with the experimental properties of electrodynamics and the weak interaction. Moreover, it contains exactly $25$ free parameters. 
    The main mathematical technique employed is to use the structure of the Clifford algebra to describe a symmetry-breaking within the 
  underlying
  matrix algebra 
   (or pair of matrix algebras).
  
   I have suggested a possible way to incorporate a unification with the strong force, but this
  is somewhat speculative, in the sense that it requires adopting the point of view that the strong force is not so much a force, as a description of the quantisation of space itself. Moreover, it requires some modification to the gauge group $SU(3)$, which may well be sufficient to rule out
  such an approach on experimental grounds.
  Further unification with general relativity is hinted at by the coincidence of two groups, both isomorphic to $SL(4,\mathbb R)$. This coincidence suggests 
  the possibility of using either the Clifford algebra $Cl(3,3)$ or perhaps a larger Clifford algebra such as $Cl(4,3)$ 
  or $Cl(4,4)$ to include models of both particle physics and gravity. I make no attempt to build such a model, and again there may be
  good physical reasons why such a model cannot exist.

\end{document}